\documentclass[12pt]{amsart}

\usepackage{amssymb}

\begin{document}
\baselineskip=15.5pt

\title[A remark]{A remark on ``Connections and Higgs fields
on a principal bundle''}

\author[I. Biswas]{Indranil Biswas}

\address{School of Mathematics, Tata Institute of Fundamental
Research, Homi Bhabha Road, Bombay 400005, India}

\email{indranil@math.tifr.res.in}

\author[C. Florentino]{Carlos Florentino}

\address{Departamento Matem\'atica, Instituto Superior
T\'ecnico, Av. Rovisco Pais, 1049-001 Lisbon, Portugal}

\email{cfloren@math.ist.utl.pt}



\date{}

\begin{abstract}
We show that a unipotent vector bundle on a non--K\"ahler
compact complex manifold does not admit a flat holomorphic
connection in general. We also construct examples of
topologically trivial stable vector bundle on compact
Gauduchon manifold that does not admit any unitary flat
connection.
\end{abstract}

\maketitle

Let $M$ be a compact connected complex manifold. A
holomorphic vector bundle $E\, \longrightarrow\, M$ is
called \textit{unipotent} if there is a filtration of
holomorphic subbundles
\begin{equation}\label{e1}
0\, =\, E_0\, \subset\, E_1\, \subset\, \cdots
\, \subset\, E_{\ell -1} \, \subset\, E_\ell \, =\, E
\end{equation}
such that for every $i\, \in\, [1\, ,\ell]$, the quotient
$E_i/E_{i-1}$ is a holomorphically trivial vector bundle.
If $M$ is K\"ahler, and $E$ is unipotent as above, then
there is a flat holomorphic connection on $E$ that preserves
each subbundle in \eqref{e1} \cite[p. 21, Corollary 1.2]{BG}
(take the Higgs field in \cite[Corollary 1.2]{BG} to be zero);
see \cite{FL} for related results.

It is natural to ask whether the same result holds for
unipotent bundles on compact complex manifolds. We will
construct an example showing that it does not hold.

Fix integers $m\, ,n\, \geq\, 1$, and fix $\tau\, \in\,
\mathbb C$ with ${\rm Im}\, \tau\, > \, 0$. This gives the
elliptic curve $T\, :=\, {\mathbb C}/({\mathbb Z}\oplus\tau
\cdot{\mathbb Z})$. Let $M$ be the corresponding
Calabi--Eckmann manifold \cite{CE}. We recall that $M$ is
diffeomorphic to $S^{2m+1}\times S^{2n+1}$, and it is the
total space of a holomorphic principal $T$--bundle over
${\mathbb C}{\mathbb P}^m\times {\mathbb C}{\mathbb P}^n$.
The complex manifold $M$ does not admit any K\"ahler metric
because $H^2(M,\, {\mathbb Z})\,=\, 0$.

Let ${\mathcal O}_M$ be the sheaf of holomorphic functions
on $M$. A computation of Borel shows that $H^1(M,\,
{\mathcal O}_M) \,=\, \mathbb C$ \cite[p. 216,
Theorem 9.5]{Bo} (see also \cite[p. 232]{Ho}). Let
\begin{equation}\label{e2}
0\,\longrightarrow\, {\mathcal 
O}_M\,\stackrel{\iota}{\longrightarrow}\, 
E \,\stackrel{p}{\longrightarrow}\, {\mathcal 
O}_M\,\longrightarrow\, 0
\end{equation}
be a nontrivial extension given by some nonzero element
of $H^1(M,\, {\mathcal O}_M)$.

Since the extension in \eqref{e2} is nontrivial, it can be
shown that the holomorphic vector bundle $E$ is nontrivial.
Indeed, if $E$ is holomorphically trivial, then take any
holomorphic section
$$
s\, :\, {\mathcal O}_M\,\longrightarrow\, E
$$
that is not a multiple of the section given by $\iota$
in \eqref{e2}. Since the composition $p\circ s$ is nonzero,
it must be an automorphism of ${\mathcal O}_M$. Hence
$s$ generates a line subbundle of $E$ which splits the short
exact sequence in \eqref{e2}. But this exact sequence does
not split. Therefore, we conclude that
$E$ is not holomorphically trivial.

Consequently, $E$ is a nontrivial unipotent vector bundle.
But $E$ does not admit any flat holomorphic connection
because $M$ being simply connected, any holomorphic
vector bundle on $M$ admitting a flat holomorphic
connection is holomorphically trivial, while we know
that $E$ is not holomorphically trivial.

Consider the short exact sequence of sheaves on $M$
$$
0\, \longrightarrow\, 2\pi\sqrt{-1}\cdot{\mathbb Z}\,
\longrightarrow\, {\mathcal O}_M\, \stackrel{\exp}{
\longrightarrow}\,{\mathcal O}^*_M \, \longrightarrow\, 0\, ,
$$
where ${\mathcal O}^*_M$ is the multiplicative sheaf of
nowhere vanishing holomorphic functions. Since
$H^i(M, \, 2\pi\sqrt{-1}\cdot{\mathbb Z})\, =\, 0$
for $i\, =\, 1\, ,2$, from the long exact sequence of
cohomologies for this short exact sequence it follows that
$$
\text{Pic}(M)\, :=\, H^1(M,\,
{\mathcal O}^*_M) \, =\, H^1(M,\,
{\mathcal O}_M) \, =\, \mathbb C\, .
$$
Since $\text{Pic}(M)$ is connected, all holomorphic
line bundles on $M$ are topologically trivial.

Therefore, any nontrivial holomorphic line bundle on $M$
is topologically trivial but does not admit any
unitary flat holomorphic connection (because $M$ being
simply connected there are no nontrivial flat
holomorphic connection on $M$).

Any holomorphic line bundle on $M$ is stable with respect
to any Gauduchon metric on $M$. Consequently, $M$ has
topologically trivial stable vector bundles that do not admit
any unitary flat holomorphic connection. (Compare this with
\cite{Bi}.)

We recall that the Hitchin--Kobayashi correspondence
implies that any stable vector bundle $E$ with $c_1(E)
\, =\, 0\, =\, c_2(E)$ on a compact K\"ahler manifold
admits a unique unitary flat holomorphic connection
(here $c_i$ are rational Chern classes).

\medskip
\noindent
\textbf{Acknowledgements.}\, The first author wishes to thank
Instituto Superior T\'ecnico, where the work was carried out,
for its hospitality. The visit to IST was funded by the FCT project
PTDC/MAT/099275/2008.

%%%%%%%%%%%%%%%%%%%%%%%%%%%%%%%%%%%%%%%%%%%%%%%%%%%%%%%%%%%%%%%%%

\end{document}